\documentclass[12pt, twoside]{article}
\usepackage{amsmath,amsthm}
\usepackage{amssymb,latexsym}
\usepackage{enumerate}
\usepackage{natbib}
\usepackage[T1]{fontenc}

\theoremstyle{definition}

\newcommand{\xx}{{\sf X}}

\newcommand{\pp}{{\sf P}}

\begin{document}

\title{Besov regularity of stochastic measures}

\author{Vadym Radchenko\\
Kyiv National Taras Shevchenko University\\
E-mail: vradchenko@univ.kiev.ua
}

\date{}

\maketitle


\renewcommand{\thefootnote}{}

\footnote{2010 \emph{Mathematics Subject Classification}: Primary 60G17; Secondary 60G57}

\footnote{\emph{Key words and phrases}: Stochastic measure, Besov spaces}



\begin{abstract}
We prove that continuous paths of $\sigma$-additive in probability set function belong to Besov space.
\end{abstract}

The Besov regularity of trajectories has been studied for some special classes of stochastic processes. Brownian motion
trajectories have been studied in~\citet{roynet}. Other Gaussian processes were considered in~\citet{cikero}. The Besov
regularity of indefinite Skorohod integral w.r.t. fractional Brownian motion was studied in~\citet{laoutu},
\citet{nuaouk}.

In the given note we consider the class of continuous stochastic processes generated by values of stochastic
measures and prove the Besov regularity of their paths.

Let $L_0=L_0(\Omega,\ \mathcal{F},\ {\pp})$ be a set of all real-valued random variables defined on the
probability space $(\Omega,\ \mathcal{F},\ {\pp})$ (more precisely, the set of equivalence classes).
Convergence in $L_0$ means the convergence in probability. Let ${\xx}$ be an arbitrary set and
${\mathcal{B}}$ be a $\sigma$-algebra of subsets of ${\xx}$.

\medskip
\textbf{Definition} {\em The $\sigma$-additive mapping $\mu:\
{\mathcal{B}}\to L_0$ is called a {\em stochastic measure}.}
\medskip

In other words, $\mu$ is a vector measure with values in $L_0$. We do not assume positivity or moment existence for
$\mu$. In~\citet{kwawoy} such $\mu$ is called a general stochastic measure.

Examples of stochastic measures are the following.

Let ${\xx}=[a,\ b]\subset\mathbb{R}_+$,
${\mathcal{B}}$ be the $\sigma$-algebra of Borel subsets of
$[a,\ b]$, and $W(t)$ be the Brownian motion. Then
$\mu(A)=\int_a^b I_A(t)\,dW(t)$ is a stochastic measure on
${\mathcal{B}}$. Moreover by the same way any continuous square
integrable martingale $X(t),\ a\le t\le b$ defines a stochastic
measure $\mu$ on ${\mathcal{B}}$ so that $\mu((s,\
t])=X(t)-X(s)$. If $W^H(t)$ is a fractional Brownian motion with
Hurst index~$H>1/2$ and $f: [0,\ T]\to\mathbb{R}$ is a bounded
measurable function then $\mu(A)=\int_0^T f(t)I_A(t)\,dW^H(t)$
is a stochastic measure on ${\mathcal{B}}$ too \citep[this fact
follows from Theorem~1.1 of][]{memiva}. Other examples may be
found in subsection 7.2 of~\citet{kwawoy}.

Let us consider an arbitrary stochastic process $X(t),\ a\le t\le b$. Put $\mu ((s,\ t])=X(t)-X(s)$ and extend $\mu$ to
algebra $\mathcal{B}_0$ of all finite unions $\cup_{k=1}^{l}\left(a_k,\ b_k\right]\subset (a,\ b]$ by additivity. Then
$\mu$ can be extended to a stochastic measure on the Borel $\sigma$-algebra iff both the following conditions holds:

(i) $\mu(A_n)\stackrel{\pp}{\to} 0$ for any $A_n\in\mathcal{B}_0,\ A_n\downarrow \emptyset$,

(ii) the set of random variables $\left\{\mu(A_n),\ n\ge 1 \right\}$ is bounded in probability for any disjoint
$A_n\in\mathcal{B}_0$
\newline \citep[Theorem~1 of][]{radtvp}.

It is known that the for any ${\xx},\ \mathcal{B}$ the set of values of any stochastic measure is bounded in
probability, i.~e.
\[
\lim_{c\to\infty}\sup_{A\in\mathcal{B}}\pp\left(\left|\mu(A)\right|>c\right)=0
\]
\citep[see][]{talagr}. Furthermore Theorem~7.1.2 of~\citet{kwawoy} establishes that the set
\begin{equation}
\label{eqbsfp} \left\{\sum_{k=1}^{n} c_{k} \mu\left(A_{k}\right),\ A_{k_1}\cap A_{k_2}=\emptyset\ \textrm{for}\ k_1\ne
k_2,\ n\ge 1,\ \left|c_{k}\right|\le 1\right\}
\end{equation}
is bounded in probability.

We consider the {\em Besov space} $B^\alpha_{pq}([a,\ b]),\ [a,\ b]\subset\mathbb{R}$. Recall that the norm in this
classical space for $1\le p,q< \infty$ and $0<\alpha< 1$ may be introduced by
\[
\|f\|^{\alpha}_{p,q}=\|f\|_{L_{p}([a,\ b])}+\left(\int_0^{b-a}\ {\left(w_p(t,\ f)\right)^q}{t^{-\alpha q-1}}
\,dt\right)^{1/q},
\]
where
\[
w_p(t,\ f)=\sup_{|h|\le t}\left(\int_{I_h} \left|f(x-h)-f(x)\right|^p\,dx\right)^{1/p},\quad
I_h=\{x\in [a,\ b]:\ x-h\in [a,\ b]\}. 
\]
 The norm in the Besov space $B^\alpha_{pp}([a,\ b])$ is equivalent to the norm in the {\em
Slobodeckij space} $W^{\alpha}_p$.

The main result of the paper will be based on Corollary~3.3 of~\citet{kamont} and the following property of stochastic
measures.

\textbf{Lemma} {\em Let $\mu$ be a stochastic measure and $a_n,\ n\ge 1$, be a sequence of positive numbers such that
$\sum_{n=1}^{\infty} a_n < \infty$. Let $\Delta_{kn}\in \mathcal{B},\ n\ge 1,\ 1\le k\le l_n$, be such that for each
$n$ and $k_1\ne k_2\ \Delta_{k_1 n}\cap \Delta_{k_2 n}=\emptyset$. Then
\[
\sum_{n=1}^{\infty} {a_n^2}\sum_{k=1}^{l_n} \mu^2\left(\Delta_{kn}\right) <\infty\quad \mbox{a.~s.}
\]
}

\textbf{Proof.} Suppose that
\[
\pp\left[\sum_{n=1}^{\infty} a_n^2 \sum_{k=1}^{l_n} \mu^2(\Delta_{kn})=+\infty\right]=\varepsilon_{0}>0.
\]
We find that for any $c>0$ there exists $j$ such that
\begin{equation}
\label{eqdefj} \pp\left[\sum_{n=1}^{j} a_n^2 \sum_{k=1}^{l_n} \mu^2(\Delta_{kn})\ge
c\right]\ge\varepsilon_{0}/2.
\end{equation}
For $j$ from~(\ref{eqdefj}), let us consider the set
\[
\Omega_1=\left\{\omega\in\Omega:\ \sum_{n=1}^{j} a_n^2 \sum_{k=1}^{l_n} \mu^2(\Delta_{kn})\ge c\right\},
\]
and the set of independent symmetric Bernoulli random variables $\varepsilon_{kn},\ 1\le n\le j,\ 1\le k\le l_n,$
defined on other probability space $(\Omega',\ \mathcal{F}',\ \pp')$, $\pp'(\varepsilon_{kn}=1)
=\pp'(\varepsilon_{kn}=-1)=1/2$. We have the following consequence of Paley-Zigmund inequality
\[
\pp'\left[\left(\sum_{1\le n\le j,\ 1\le k\le l_n} \lambda_{kn}
\varepsilon_{kn}\right)^2\ge({1}/{4})\sum_{1\le n\le j,\ 1\le k\le l_n} \lambda_{kn}^2 \right]\ge {1}/{8},\quad
\lambda_{kn}\in \mathbb{R}
\]
(see, for example, Lemma~V.4.3~(a)~\citet{vatach} or Lemma~0.2.1~\citet{kwawoy} for~$\lambda=1/4$).

By applying this inequality with taking $\lambda_{kn}=a_n \mu(\Delta_{kn},\ \omega)$ for a fixed $\omega\in\Omega_1$ we
get
\[
\pp'\left[\omega':\ \left(\sum_{n=1}^{j} a_n \sum_{k=1}^{l_n} \varepsilon_{kn}(\omega')\mu(\Delta_{kn},\
\omega)\right)^2\ge {c}/{4}\right]\ge{1}/{8}.
\]
Integrating the above inequality with respect to measure $\pp$ over $\Omega_1$ we obtain
\[
\pp\times\pp'\left[(\omega,\ \omega'):\ \left(\sum_{n=1}^{j} a_n \sum_{k=1}^{l_n}
\varepsilon_{kn}(\omega')\mu(\Delta_{kn},\ \omega)\right)^2\ge{c}/{4}\right]\ge {\varepsilon_{0}}/{16}.
\]
Hence by using Fubini's theorem, we have that there exists $\omega'_0\in\Omega'$ such that
\[
\pp\left[\omega:\ \left(\sum_{n=1}^{j} a_n \sum_{k=1}^{l_n} \varepsilon_{kn}(\omega'_0)\mu(\Delta_{kn},\
\omega)\right)^2\ge {c}/{4}\right]\ge{\varepsilon_{0}}/{16}.
\]
Recalling that each $\varepsilon_{kn}(\omega'_0)=1$ or $\varepsilon_{kn}(\omega'_0)=-1$ we obtain sets $B_n,\
C_n\in\mathcal{B}$ such that
\begin{equation}
\label{eqsuub} \pp\left[\left|\sum_{n=1}^{j} a_n (\mu(B_{n})-\mu(C_{n}))\right|\ge \sqrt{c}/{2}\right]\ge
{\varepsilon_{0}}/{16}.
\end{equation}
We have
\begin{equation}
\label{eqsubo}
 \max_{x\in\xx}\left|\sum_{n=1}^{j} a_n\left(1_{B_{n}}(x)-1_{C_{n}}(x)\right)\right|\le \sum_{n=1}^{\infty} a_n.
\end{equation}
Recall that $\varepsilon_{0}>0$ is fixed and $c>0$ is arbitrary. Therefore (\ref{eqsuub}) and (\ref{eqsubo}) contradict
the boundedness in probability of the sums~(\ref{eqbsfp}). This completes the proof of the Lemma. \hfill$\square$

The main result of the paper is the following.

\textbf{Theorem} {\em Let $\xx=[a,\ b]\subset\mathbb{R}$, $\mathcal{B}$ be the Borel $\sigma$-algebra, $\mu$
be a stochastic measure on $\mathcal{B}$ and the process $\mu(t)=\mu([a,\ t]),\ a\le t\le b$, have continuous paths.
Then for any $p\ge 2,\ 0<\alpha< 1/p$, the path of $\mu(t)$ with probability~1 belongs to the Besov space
$B^\alpha_{pp}([a,\ b])$.}

\textbf{Proof.} When $f:[a,\ b]\to\mathbb{R}$ is a continuous function, Corollary~3.3 of~\citet{kamont} shows that the
convergence of the series
\[
\sum_{n=1}^{\infty}2^{ n (\alpha p-1)}\sum_{k=1}^{2^n}\left|f\left(a+k2^{-n}(b-a)\right)-
f\left(a+(k-1)2^{-n}(b-a)\right) \right|^p
\]
implies that $f\in B^\alpha_{pp}([a,\ b])$. Obviously, it is sufficient to prove the convergence for the second power
of differences of function values.

By taking
\[
a_n=2^{ n (\alpha p-1)/2},\quad \Delta_{kn}=\left(a+(k-1)2^{-n}(b-a),\ a+k2^{-n}(b-a)\right],\quad
  1\le k\le 2^n
\]
in the Lemma, we see that continuous paths of~$\mu(t)$ a.~s. satisfy the mentioned
condition.

\hfill$\square$

In particular the statement of the Theorem may by applied to the paths of any continuous square integrable martingale.

Embeddings of the Besov spaces \citep[see, for example, subsection 3.2.4 of][]{triebel} implies that a continuous path
of $\mu(t)$ a.~s. belongs to $B^\alpha_{pq}([a,\ b]),\ q\ge p\ge 2,\ 0<\alpha< 1/p$.

Note that with probability~1 the trajectories of the Brownian motion do not belong to the Besov spaces
$B^\alpha_{pq}([0,\ 1])$ for all $1/2<\alpha<1,\ p,\ q\ge 1$ \citep[Theorem~1 of][]{roynet}.

\section*{Acknowledgments}
This work was supported by Alexander von Humboldt Foundation, grant 1074615.
I am grateful to Prof. M.~Z\"{a}hle-Ziezold for fruitful discussions during the preparation of this paper.

\end{document}